\documentclass[12pt]{article}
\newtheorem{ttt}{Theorem}
\newtheorem{ddd}{Definition}[section]

\newtheorem{lll}{Lemma}[section]

\newcommand{\Th}{$\Theta$}
\newcommand{\aut}{automorphism\hspace{2pt}}

\title{The group of automorphisms of the category of free associative algebras}
\author{A. Berzins\\University of Latvia\\ e-mail: aberzins@latnet.lv}
\date{August 2004}

\begin{document}
\maketitle



\begin{abstract}

In this paper, the problem formulated in [8] is solved. We prove,
that the group of automorphisms of the category of free
associative algebras is generated by semi-inner and mirror
automorphisms.

\end{abstract}

\section{Introduction}

The classical algebraic geometry arise from consideration systems
of polynomial equations in an affine space $P^n$, where $P$ is
arbitrary field. Although equations solving in algebras in
different varieties of algebras $\Theta$ was considered long ago,
foundation of the universal algebraic geometry are B. Plotkin's
scientific works [5-8] where concepts of affine space, algebraic
set, algebraic variety for arbitrary varieties of algebras were
introduced.

The basis of this theory is consideration of a point as an element
from ${\rm Hom}(W,X)$, where $W$ is a free finitely generated
algebra, and $X$ is a fixed algebra in variety of universal
algebras \Th. It was found that many problems in this theory, such
as geometric equivalence, geometric similarity, isomorphism and
equivalence of categories of algebraic sets and varieties depend
on structure of ${\rm Aut}(\Theta^0)$ and Aut(End($W$)), where
$\Theta^0$ is the category of algebras $W=W(X)$ with the finite
$X$ that are free in $\Theta$, and $W$ is a free algebra in
$\Theta$. Structure of ${\rm Aut(Com}-P))^0$ (i.e, the classical
case) was described in [1]. Also are described ${\rm
Aut}(\Theta^0)$ for varieties of groups, semigroups, Lie algebras
and some another varieties. For some varieties are described
Aut(End($W$)), but these problems are more difficult. In
particular, this problem is not solved even in the classical case
of free commutative algebras, i.e. the structure of
Aut(End$(P[x_1,\ldots ,x_n])$) is not describe if $n>2$. In all
solved cases groups Aut(End($W$)) are generated by semi-inner
automorphisms and mirror automorphism.

In this article we describe ${\rm Aut(Ass}-P)^0$. When the paper
had been finished, I was informed that a solution of this problem
is contained in the forthcoming paper [4]. Both solutions of this
problem are principally different, and methods of proves are
independent. Explicitly I want to note that proof in the present
paper does not use the description of the group ${\rm
Aut}(W(x,y)$). It is very important, because structure of the
group ${\rm Aut}(W(x_1,\ldots ,x_n))$ is not described, when
$n>2$, but counterexamples show that structures of ${\rm
Aut}(W(x_1,\ldots ,x_n))$ and ${\rm Aut}(W(x,y)$) are principally
different. So, methods presented in the paper may be useful to
describe the group Aut(End($W_n$)) for $n>2$.

\section{Definitions}

Let $(\textrm{Ass}-P)^0=\Theta^0$ be the category of free
associative algebras over a field $P$ with a finite set of
generators, and $\tau\in \textrm{Aut}(\Theta^0)$. Then there exist
a set of bijections $\mu=\{\mu_i\}$, $i\in \textrm{N}$ of $W_i$,
such that for every $s\in \textrm{Hom}(W_i,W_j)$
$$s^\tau=\mu_i s \mu_j^{-1}$$
Here $W_k$ is the free \textit{k}-generated associative algebra.
Automorphism $\tau$ which is representable in the form
$s^{\tau}=\mu s\mu^{-1}$ is called quasi-inner. In the category
$(\textrm{Ass}-P)^0$ every automorphism is quasi-inner.  Of
course, arbitrary set of substitutions does not generate an
automorphism. For category $\Theta^0$ we have three types of
bijections $\mu$ which generate an automorphism.\\

1. $\mu=\overline{\alpha}$, where $\overline{\alpha}$ is the
natural extension of an automorphism $\alpha \in \textrm{Aut}(P)$
on $W$, i.e.
$\overline{\alpha}\left( \sum a_ku_k\right)=\sum \alpha(a_k)u_k$.
 We shall write $\overline{\alpha}=\alpha$ and call it "\aut of
 field".Here $u_k$ is a product of variables.\\

 2. $\mu=\eta=\{\eta_i\}$, $\eta_i\in \textrm{Aut}(W_i)$.\\

Let as recall some definitions for Ass-\textit{P}. (See [8] for
general case)

\begin{ddd}

The automorphism $\tau$ generated by a set of automorphisms $\eta
= \{\eta_i\in W_i\}$ is called inner \aut.

\end{ddd}

\begin{ddd}

Note, that every $\alpha \in \textrm{Aut}(P)$ belongs to
normalizer of subgroup  $\textrm{Aut}(W_i)$ in the group of all
bijections of $W_i$. So, every product of elements of
$\textrm{Aut}(P)$ and $\textrm{Aut}(W_i)$ my be represented in the
form $\mu_i=\alpha\eta_i$, $\alpha \in \textrm{Aut}(P)$,
$\eta_i\in\textrm{Aut}(W_i)$. Automorphism $\tau \in
\textrm{Aut}(\Theta^0)$ generated by $\mu$ is called semi-inner.

\end{ddd}

3. Now we describe the third type of $\mu$. Let $S=S(X)$ be a free
semi-group. For every $u=x_{i_1}x_{i_2}\ldots x_{i_n}$ in $S$,
take $\overline{u}=x_{i_n}\ldots x_{i_2}x_{i_1}$. Then
$u\rightarrow \overline{u}$ is an anti\aut of the semigroup S.

\begin{ddd}

Let $W=W(X)$ be a free associative algebra. For every its element
$\omega=\lambda_0+\lambda_1u_1+\cdots+\lambda_ku_k$ denote
$\overline{\omega}=\lambda_0+\lambda_1\overline{u_1}+\cdots+\lambda_k\overline{u_k}$.
The transition $\omega\rightarrow\overline{\omega}$ is an
antiautomorphism of the algebra $W$. The set of these transitions
$\beta=\{\beta_i:W_i\rightarrow W_i\}$ generate the automorphism
$\delta\in \textrm{Aut}(\Theta^0)$, $s^\delta =\beta s
\beta^{-1}$. We call $\beta$ mirror antiautomorphism of $W$ and
$\delta$ mirror automorphism of $\textrm{Aut}(\Theta^0)$.

\end{ddd}

Here $\delta$ is not inner and is not semi-inner, but is
quasi-inner. Note that

\begin{enumerate}

\item $\delta$ belongs to the normalizer of subgroup
$\textrm{Sinn}(\Theta^0)$ in $\textrm{Aut}(\Theta^0)$ that
consists of all semi-inner automorphisms;

\item every antiautomorphism of $W$ is a product of the mirror
antiautomorphism $\beta$ and an automorphism of $W$.

\end{enumerate}

\section{Almost central algebras}

Let $W(X)$ be a free finitely generated algebra in variety
$\Theta$.

\begin{ddd}

Bijection $\mu : W(X)\mapsto W(X)$ is called central, if for every
$s\in \textrm{End}(W(X))$
$$\mu s=s\mu.$$
Algebra $W(X)$ is called central, if every its central bijection
is identical.

\end{ddd}

\begin{ttt}

Let $W(X)$ be a free central algebra and $\mu$ some bijection of
$W(X)$ generating an automorphism
$\tau\in\textrm{Aut}(\textrm{End}(W(X)))$. Then $\mu$ transforms
every base of algebra $W(X)$ to a base of algebra.

\end{ttt}
\textbf{Proof.} [see 3]\\[6pt]

For free algebras in a variety $\Theta$ of algebras over field $P$
(Com-P, Ass-P, and others) we shall correct the definition of
central algebra. Clearly, in every such algebra $W(X)$ we have a
bijection $\mu(u)=au+b$, where $a,b\in P,$ $a\neq 0$ (linear
bijection). This bijection commutate with every endomorphism of
algebra $W(X)$ and, obviously, transforms every base of $W(X)$ to
a base.

\begin{ddd}

Finitely generated algebra $W(X)$ over field $P$ is called almost
central, if every its central bijection is linear.

\end{ddd}

\begin{ttt}

Free finitely generated commutative and associative algebras are
almost central.

\end{ttt}
\textbf{Proof.} Let $\mu$ be a central bijection of $W(X)$,
$\mu(x_i)=r_i(x_1,x_2,\ldots,x_n)\in W(X)$ , and $s\in
\textrm{End}(W(X)$, $s(x_1)=x_1$, $s(x_i)=0$ for $i=1$. Then
$$\mu s(x_1)=\mu (x_1)=r_1(x_1,x_2,\ldots,x_n),$$
$$s\mu(x_1)=s(r_1(x_1,x_2,\ldots,x_n))=r_1(x_1,0,\ldots,0)=r(x_1).$$
So $r_1(x_1,x_2,\ldots,x_n)=r(x_1)$ is a polynomial of one
variable.

For arbitrary $u\in W(X)$ take $s\in \textrm{End}(W(X))$
$s(x_1)=u$, $s(x_i)=x_i$ for $i>1$. Then
$$\mu(u)=\mu s(x_1)=s\mu(x_1)=s(r(x_1))=r(u),$$
and so $\mu(u)=r(u)=a_0+a_1+\cdots+a_ku^k$. Because $\mu$ is
bijection, so $r$ is linear, Q.E.D..\\[3pt]
\textbf{Corollary.} From the theorems 1 and 2 we have that in
mentioned algebras every bijection $\mu$, which generate an
automorphism of $\textrm{Aut}(\textrm{End}(W(X)))$  transforms
every base of algebra $W(X)$ to a base of algebra.

\section{The group ${\rm Aut(Ass}-P)^0$}

Now we can formulate the main result of this article (problem 6 in
[8]).

\begin{ttt}

The group $\textrm{Aut}(\textrm{Ass}-P)^0$ is generated by
semi-inner and mirror automorphisms.

\end{ttt}
\textbf{Proof.} Let $\tau$ be an automorphism of $\Theta^0$. Then
exist a set of bijections $\mu={\mu_i}$, such that for every $s\in
\textrm{Hom}(W_i,W_j)$ equality $s^{\tau}=\mu_js\mu_i^{-1}$ holds.
We must prove that $\mu_i=\alpha\cdot\eta_i$ for every $i\in
\textrm{N}$ or $\mu_i=\alpha\cdot\eta_i\cdot\beta$ for every
$i\in \textrm{N}$. Here $\alpha\in \textrm{Aut}(P)$, $\eta_i\in
\textrm{Aut}(W_i)$ and $\beta$ is mirror.

Let us consider the free 1-generated associative algebra
$W(x)=P[x]$, i.e. the free 1-generated commutative algebra. From
the analogous theorem for $\textrm{Aut}(\textrm{Com}-P)^0$ (see
[1]) we have that $\mu_1=\alpha\eta_1$. So, we can consider
$\mu'=\alpha^{-1}\mu$. Then $\mu'_1$ is automorphism of $W(x)$.
Clearly, the theorem will be proved, if we check that $\mu'$ is a
set of automorphisms or a set of antiautomorphisms of $W_i$.

Write $\mu=\mu'$. Let
$s\in\textrm{Hom}(W_1,W_2),\hspace{5pt}s:t\rightarrow f(x,y)$. It
follows from $\mu_2s=s^{\tau}\mu_1$ that $\mu_2$ is a morphism of
1-generated subalgebra $P(f(x,x))$ of $W(x,y)$ to $W(x,y)$.

Let $s\in S=\textrm{Hom}(W(x,y),W(t))$. Then
$\mu_1s=s^{\tau}\mu_2$; $\mu_1s(x)=s^{\tau}\mu_2(x)$,
$\mu_1s(y)=s^{\tau}\mu_2(y)$, $\mu_1s(x+y)=s^{\tau}\mu_2(x+y)$.
Denote $\mu=\mu_2$. Then
$$s^{\tau}(\mu(x+y)-\mu(x)-\mu(y)=\mu_1s((x+y)-x-y)=0,$$
and
\begin{equation}
\mu(x+y)-\mu(x)-\mu(y)\in \bigcap_{s\in S}\textrm{Ker}s=(xy-yx)
\end{equation}
Here $(xy-yx)$ is ideal generated by $xy-yx$. Analogously

\begin{equation}
\mu(x\cdot y)-\mu(x)\cdot\mu(y)\in (xy-yx)
\end{equation}

\begin{lll}
$\mu(u+v)=\mu(u)+\mu(v)$ for every $u,v\in W(X)$.
\end{lll}
\textbf{Proof.} Let
$$\Delta(x,y)=\mu(x+y)-\mu(x)-\mu(y)=f(x,y)=g(\mu(x),\mu(y).$$
Last equality holds because by theorem 2 algebra $W(X)$ is almost
central. Note that $g(u,v)$ is homogenous of degree 1. Really,
$$g(\alpha\mu(x),\alpha\mu(y))=g(\mu(\alpha x),\mu(\alpha y))=$$
$$\mu(\alpha x+\alpha y)-\mu(\alpha x) -\mu(\alpha y)= \alpha\cdot
g(\mu(x),\mu(y)).$$ So $\mu(x+y)-\mu(x)-\mu(y)=a\mu(x)+b\mu(y)$,
i.e. $\mu(x+y)=c\mu(x)+d\mu(y)$. Using endomorphism we have
$\mu(u+v)=c\mu(u)+d\mu(v)$ for every $u,v\in W(X)$. Substitutions
$u=0$ and $v=0$ gives that $c=d=1$, Q.E.D.

\begin{lll}
$\mu(xy)=\mu(x)\mu(y) +a(xy-yx)$, where $a\in P$.
\end{lll}
\textbf{Proof.} Let $\Delta(x,y)=\mu(xy)-\mu(x)\cdot\mu(y).$ Then
$\Delta(x,y)$ is $x$ and $y$ homogenous of degree 1. Since
$\Delta(x,y)\in(xy-yx)$ (formula (2)), so $\Delta(x,y)=a(xy-yx)$.
Q.E.D.
\begin{lll}
$\mu(xy)=\mu(y)\mu(x) +b(xy-yx)$, where $b\in P$.
\end{lll}
The proof is the same as the proof of the lemma 4.2.
\begin{lll} $\mu(xy)=\alpha\mu(x)\mu(y)+\beta\mu(y)\mu(x)$, where
$\alpha,\beta\in P$.
\end{lll}
\textbf{Proof.} It follows from lemma 4.2 and lemma 4.3 that
$$(b-a)\mu(xy)=b\mu(x)\mu(y)-a\mu(y)\mu(x)$$
If $a=b$, then from lemma 4.2 and lemma 4.3 we get equality
$\mu(x)\mu(y)=\mu(y)\mu(x)$ which easy leads to contradiction. So
we can take $\alpha=\frac{b}{b-a}$ and $\beta=\frac{-a}{b-a}$ and
get required equality.\\[6pt]

\noindent\textbf{Corollary.} For every $u,v\in W(x,y)$
\begin{equation}
\mu(uv)=\alpha\mu(u)\mu(v)+\beta\mu(v)\mu(u)
\end{equation}
\begin{lll}
$\mu(uv)=\mu(u)\mu(v)$ for every $u,v\in W(x,y)$ or\\[3pt]
$\mu(uv)=\mu(v)\mu(u)$ for every $u,v\in W(x,y)$.
\end{lll}

\noindent\textbf{Proof.} Let $u,v$ be elements of $W(x,y)$, such
that $\mu(u)=x$ and $\mu(v)=y$. Let us denote
$\mu(r)=\overline{r}$. Then
$$\mu(uuv)=\alpha\overline{uu}\cdot\overline{v}+\beta\overline{v}\cdot\overline{uu}=\alpha
xxy+\beta yxx,$$
$$\mu(uuv)=\alpha\overline{u}\cdot\overline{uv}+\beta\overline{uv}\cdot\overline{u}=\alpha^2xxy+2\alpha\beta
xyx+\beta^2yxx.$$ Therefore $\alpha^2=\alpha$, $\alpha\beta=0$,
$\beta^2=\beta$. Because the case $\alpha=\beta=0$ is impossible,
we have two solutions:
\begin{enumerate}
\item
$\beta=0,\alpha=1$, then $\mu(uv)=\mu(u)\mu(v)$;
\item
$\beta=1,\alpha=0$, then $\mu(uv)=\mu(v)\mu(u)$.
\end{enumerate}

It follows from $\mu(\alpha u)=\alpha\mu(u)$, lemma 4.1 and lemma
4.5 that $\mu$ is an automorphism or antiautomorphism of $W(x,y)$.
Using endomorphisms $s\in\textrm{Hom}(W_2,W_n)$ it is easy to
check that $\mu_n$ also is an automorphism or antiautomorphism.
The theorem is proved.\\[12pt]
{\Large\textbf{Acknowledgments}}\\

The author is happy to thank professor B. Plotkin for stimulating
discussions of the results.\\[18pt]
{\Large\textbf{References}}
\begin{enumerate}

\item A. Berzins, {\it The automorphisms of \,{\rm End}{}$K[x]$}, Proc.
Latvian Acad. Sci., Section B, 2003, vol. 57, no. 3/4, pp. 78-81.

\item A. Berzins, {\it Geometric equivalence of algebras}, Int. J.
Alg. Comput., 2001, vol. 11, no. 4, pp. 447-456.

\item A. Berzins, B. Plotkin and E. Plotkin, {\it Algebraic
geometry in Varieties with the Given Algebra of constants}, J.
Math. Sci., New York, 2000, vol. 102, no. 3, pp. 4039-4070.

\item G. Mashevitzky, {\it Automorphisms of categories of free
 associative algebras}, Preprint.

\item B. Plotkin, {\it Algebraic logic, varieties of algebras and
algebraic varieties}, Proc. Int. Alg. Conf., St. Petersburg, 1995,
Walter de Gruyter, New York, London, 1996.

\item B. Plotkin, {\it Varieties of algebras and algebraic
varieties}, Israel J. of Mathematics, 1996, vol. 96, no. 2, pp.
511-522.

\item B. Plotkin, {\it Varieties of algebras and algebraic
varieties. Categories of algebraic varieties}, Sib. Adv. Math.,
1997, vol. 7, no. 2, pp. 64-97.

\item B. Plotkin, {\it Algebras with the same (algebraic)
geometry}, Proc. Steklov Inst. Math, Vol. 242, 2003, pp. 165-196.

\end{enumerate}

\end{document}